\documentclass[3p,times]{elsarticle}

\usepackage{caption,ecrc,epstopdf}

\volume{00}

\firstpage{1}

\journalname{Physica A}

\runauth{M.~N.~\c{C}ankaya and J.~Korbel}

\jnltitlelogo{Physica A}

\usepackage{amsmath}
\usepackage{amssymb}
\usepackage{amsthm}
\usepackage{dsfont}
\usepackage{graphicx,color}
\usepackage{subcaption}
\newcommand{\ud}{\mathrm{d}}
\dochead{}
\begin{document}
\begin{frontmatter}
%Generalized Entropies and
 %%%%%%%%%%%%%%%%%%%%%%%%%%%%%%%%%%%%%%%%%%%%%%%%%%%%%%%%%%%%%%%%%%%%%%%%%%%%%%%%%%%%%%%%%%%%%%%%%%
\title{Least informative distributions in Maximum $q$-log-likelihood estimation}% Force line breaks with \\
%\title{On Statistical Properties of Jizba-Arimitsu Hybrid Entropy}% Force line breaks with \\
%%%%%%%%%%%%%%%%%%%%%%%%%%%%%%%%%%%%%%%%%%%%%%%%%%%%%%%%%%%%%%%%%%%%%%%%%%%%%%%%%%%%%%%%%%%%%%%%%%%%%
%
\author[USK,USAK]{Mehmet Niyazi \c{C}ankaya}
\ead{mehmet.cankaya@usak.edu.tr}
\author[MUW,CSH,CTU]{Jan Korbel}%
\ead{korbeja2@fjfi.cvut.cz}

\address[USK]{Applied Sciences School, Department of International Trading, U\c{s}ak, Turkey}
\address[USAK]{U\c{s}ak University, Faculty of Art and Sciences, Department of Statistics, U\c{s}ak, Turkey}
\address[MUW]{Section for the Science of Complex Systems, CeMSIIS, Medical
University of Vienna, Spitalgasse 23, 1090 Vienna, Austria}
\address[CSH]{Complexity Science Hub Vienna, Josefst\"{a}dterstrasse 39, 1080 Vienna, Austria}
\address[CTU]{Faculty of Nuclear Sciences and Physical Engineering, Czech Technical University in Prague,
B\v{r}ehov\'{a} 7, 115 19 Praha 1, Czech Republic}

%\textcolor{red}{ }
\begin{abstract}
We use the Maximum $q$-log-likelihood estimation for Least informative distributions (LID) in order to estimate the parameters in probability density functions (PDFs) efficiently and robustly when data include outlier(s). LIDs are derived by using convex combinations of two PDFs, $f_\epsilon=(1-\epsilon)f_0+\epsilon f_1$. A convex combination of two PDFs is considered as a contamination $f_1$ as outlier(s) to underlying $f_0$ distributions and $f_\epsilon$ is a contaminated distribution. The optimal criterion is obtained by minimizing the change of Maximum q-log-likelihood function when the data have slightly more contamination. In this paper, we make a comparison among ordinary Maximum likelihood, Maximum q-likelihood  estimations, LIDs based on $\log_q$ and Huber M-estimation. Akaike and Bayesian information criterions (AIC and BIC) based on $\log_q$ and LID are proposed to assess the fitting performance of functions. Real data sets are applied to test the fitting performance of estimating functions 
that include shape, scale and location parameters.
\end{abstract}
%\textbf{Highlights}%of  from the generalized function and assessing novel AIC and BIC for LIDs and PDFs

%\begin{itemize}when  we take classical derivative from a special case of variational calculus with respect to $\epsilon$ and $\epsilon=0$ for generalized functions  $\log_q(f_\epsilon)$.estimation
 % \item  These generalized functions and their LIDsand test them for fitting performanceT the contamination rate $\epsilon$  is eliminated.The elimination of $\epsilon$ is performed%PDFs in $\log_q$  are also proposed to compare with LIDs in $\log_q$.
  %\item
%\end{itemize}
\end{frontmatter}

\section{Introduction}

Least informative distributions (LIDs) under various characterizing restrictions on Fisher information were considered by
\cite{Nianotherlook,Hub64,Hub81,Shev08redes}. LIDs in estimating functions from M-estimation are  proposed in \cite{Hub64}. Forerunners of estimating functions and estimating equations can be found e.g. in Refs. \cite{God60,GodTh84}. Some more examples of estimating functions are from Refs. \cite{Nianotherlook,Hub64,Hub81,Shev08redes,LeeShevLID,SteBoos}. Maximum $q$-likelihood estimation method in logarithm with $q-$difference operator ($\log_q$) as a generalized logarithm \cite{MLqEBio,FerrariYang10,Berchermt12} can be given as an example for estimating functions. LIDs provide a special type of M-estimation, which minimizes the change of the Maximum likelihood function under increasing contamination.
 %for contaminated distribution when one wants to fit data via
%a convex combination of two probability density functions (PDFs) instead of using only one PDF in the sake of having more information from data. via  $\log_q$ \textcolor[rgb]{1.00,0.00,0.00}{}
Thus,  M-estimators from LIDs based on $\log_q$ can be not only robust but also efficient. An advantage of LID is that a neighborhood of a probability density function (PDF)  can be obtained. PDFs  have been proposed from  maximum entropy principle by using  generalized entropies \cite{CanKor16,Korbelrescaling,JKZT,KanMEPpow,MatsuHen}. The comprehensive survey on generalized divergence and their applications are introduced by \cite{cicalpbet,amaritablediv,egucetalent,tomaentropy14,kanamor,kanamor2,Eguchikat,Yalcinetal,Suyariq2,Suyariq}. A mixed PDF with fixed mixing proportions (or contamination rates) $\epsilon_1$, $\epsilon_2$ and $\epsilon_3$ is used to construct a bimodal from two mixings \cite{Hadbimod11}  and a trimodal from three mixings \cite{Hadtrimod} distributions for Ising model.  However, LID is free from the mixing proportion $\epsilon$. It is not essential to know the mixing proportions or we are not interested in the estimations of mixing proportions due to problems on estimating the mixing proportions.

M-estimation that is generalization of Maximum likelihood estimation (MLE) is used to produce robust estimators for parameters of a PDF.  M-estimators  \cite{Hub64,God60} are defined through an estimating functions  minimizing $ \sum_{i=1}^n \Lambda[f_0(x_i,\boldsymbol \theta)]$ over $\boldsymbol \theta$ \cite{Hub64,Hub81,Andrewsetalloc72,Hampeletal86}. Here,  $\Lambda$ is a concave function that is capable of making an one to one mapping from $f_0(x,\boldsymbol \theta)$ to $\Lambda$. M-estimators are derived by fixed functions, such as Huber, Tukey, etc. LIDs occurred by restriction on Fisher information are used to produce Huber, Tukey, etc \cite{Hub64}. MLE as a special case of M-estimation is a method for estimations of parameters in a PDF. It is based on logarithm and does not work properly to estimate parameters in a PDF efficiently and robustly when data which include outlier(s) are non identically distributed, therefore we will use function $\log_q$ that mimics MLE method \cite{MLqEBio,FerrariYang10}. In our proposition, the benefit of LIDs and a PDF in $\Lambda$  is that one can propose the estimating functions from arbitrary PDFs to get more precise estimators for parameters in PDFs. The more precision can be accomplished by the parameter $q$ and also LID in $\log_q$. We propose to use LIDs and a PDF in $\log_q$  to get new estimating functions and compare them with Huber M-estimation.  Finally, we have estimating functions  to fit data and information criterions for these functions  by using $\log_q$.

The rest of paper is organized as follows. Section \ref{IntroPre} is a composition of preliminaries about estimation.
Section \ref{LIDGenEntropy} proposes a PDF and LIDs in generalized logarithms as new functions in M-estimation.  Section \ref{modAIC} is provided to assess the fitting competence of estimating functions.  Section \ref{realsection} is considered for fitting competence of proposed LIDs and assessing novel Akaike and Bayesian information criterions  (AIC and BIC) for LIDs and a PDF in the function $\Lambda$. Finally, a conclusion is given in section \ref{conclusionsection}.

\section{Preliminaries}\label{IntroPre}

\subsection{Estimation procedure}

An estimation procedure is performed when one has a sampled version of a PDF, i.e, $f_0(x;\hat{\boldsymbol \theta})$, of $f_0(x;\boldsymbol \theta)$. The Maximum likelihood estimation (MLE) defined as %already uses a function in itself, as given by following form:$\mathbf{x}=\{x_1,x_2,...,x_n\}$,
\begin{equation}\label{MLE}
  L(f_0) = \prod_{i=1}^{n} f_0(x_i;{\boldsymbol \theta}),
\end{equation}
\noindent where  $\boldsymbol \theta$ is a vector for parameters, $f_0: \mathbb{R} \times \mathbb{R}^d \rightarrow [0,1]$, $\boldsymbol \theta \in \mathbb{R}^d$, $d \in \mathbb{Z}^{+}$, $n$ is a number for sample size of data set generated randomly from $f_0$ and $f_0(x;\boldsymbol \theta)$ is differentiable with respect to (w.r.t) parameters $\boldsymbol \theta$. For convenience, MLE can be also expressed as $\sum_{i=1}^{n} \log[f_0(x_i;{\boldsymbol \theta})]$  \cite{gutgraduate}. %is maximized according to parameters \cite{gutgraduate}.

M-estimation is a generalization of MLE proposed in \cite{Hub64,Hub81,God60}. An estimating function from M-estimation can be defined by following form for a PDF
\begin{equation}\label{estimatingPDF}
\rho_\Lambda(f_0)=  \sum_{i=1}^{n} \Lambda[f_0(x_i;\boldsymbol \theta )].
\end{equation}
Optimization (maximization or minus minimization) of $\rho_\Lambda$ over parameters $\boldsymbol \theta$ produces M-estimators $\hat{\boldsymbol \theta}$. If $\Lambda$ is substituted by $\log$, then $\rho_\Lambda$ becomes the ordinary MLE  \cite{Hub64,Hub81,God60,Hampeletal86}. When $\Lambda$ is $\log_q$, MLE based on $\log_q$ is obtained and it is called as the Maximum q-log-likelihood estimation (MqLE) method  \cite{MLqEBio,FerrariYang10}.

Huber M-estimation leads to set of equations $\sum_{i=1}^n \psi(y_i) = 0$, where $\psi(f_0) = \nabla_\theta \rho(f_0)$. As a typical example is the estimation of mean $\mu$ and standard deviation $\sigma$. For this end, one can choose $\psi(y)=  y $ for $ |y| \leq u$ and $\text{sign}(y)u$ for  $|y| > u$. $y=\frac{x-\mu}{\sigma}$ is a score function of normal distribution when $|y| \leq u$. $\text{sign}(y)u$ is a score function of Laplace distribution when $|y| > u$. $u$ is a tuning parameter to have robust estimators \cite{Hub64,Hub81}.

\subsection{q-deformed logarithms and connection to Tsallis entropy}
The log-likelihood estimation is based on maximization of sum of $\log f(x_i;\boldsymbol \theta)$. On the other hand, when we want to focus more on rare events, it can be convenient to generalize the log-likelihood to a special kind of M-estimator, based on generalization of logarithm. The generalized logarithm is defined as
\begin{equation}
\log_q(t)=\frac{t^{1-q}-1}{1-q},
\end{equation}
\noindent for $q > 0$. For $q \rightarrow 1$, we recover the ordinary logarithm. The aim of the generalization of MLE method is to focus on the large or small probabilities, because with changing $q$, the importance of probabilities is changed. The $q$-log-likelihood function can be therefore defined as
\begin{equation}
\rho_{\log_q}(f_0(x_i;\boldsymbol \theta)) = \sum_{i=1}^n \log_q[f_0(x_i;\boldsymbol \theta )].
\end{equation}

Interestingly, the generalized logarithm is closely related to Tsallis entropy. Tsallis entropy is a non-additive generalization of Shannon entropy $H[f_0(x_i;\boldsymbol \theta)] = -\sum_{i=1}^n f_0(x_i;\boldsymbol \theta)  \log[f_0(x_i;\boldsymbol \theta)] = E[ \log(1/f_0)]$. It is defined as
\begin{equation}
S_q(f_0(x_i;\boldsymbol \theta)) = E[ \log_q(1/f_0) ] = \frac{1}{1-q} \left(\sum_i [f_0(x_i;\boldsymbol \theta)]^q -1 \right) \, .
\end{equation}
One can find the relation between Maximum $q$-log-likelihood function and Tsallis entropy, which is
\begin{equation}
\rho_{\log_q}[f_0(x_i;\boldsymbol \theta)] = S_q[f_0(x_i;\boldsymbol \theta)^{(1-q)/q}]\, .
\end{equation}

\section{Least informative distributions based on generalized logarithms}\label{LIDGenEntropy}
Generalized entropies and connected $q$-deformed algebra have found many applications in physics and related fields. \cite{tsallis1988possible,Jizba3,JizbaKorbel16,MachadoEntropy,UbriacoFC,Wadatwopara}.

%Our aim is to use Maximum $q$-log-likelihood distribution in order to find least informative distributions.
Our aim is to propose the LIDs based on  the Maximum q-log-likelihood. Let us consider a convex combination $f_\epsilon$

\begin{equation}\label{fconvexLID}
f_\epsilon(x;\boldsymbol \theta)=(1-\epsilon)f_0(x;\boldsymbol \theta)+\epsilon f_1(x;\boldsymbol \theta)
\end{equation}
\noindent composed of PDFs $f_0$ and $f_1$, and $\epsilon \in [0,1]$, which a contamination rate. $f_\epsilon$ is a contaminated distribution with contamination $f_1$ as outlier(s) and $f_0$ is an underlying distribution. It can also considered that it represents the situation when we have non-identically distributed data.

%$q$-log-likelihood
Our aim is to find the optimal parameters $\boldsymbol \theta$ in the function $f_\epsilon$, for which the function $\Lambda$ exhibits a minimal change w.r.t the parameter $\epsilon$ when we contaminate the underlying distribution $f_0$ by a small amount of outlier distribution $f_1$, i.e., we set $\epsilon$ in $f_\epsilon$ as a small value close to zero. % From the mathematical point of view, let us define the maximum $q$-log-likelihood estimator for $f_\epsilon$
%
%\begin{equation}\label{genvariational.calculusin.Fisher}
%\rho_{\log_q}(f_0,f_1)= \sum_{i=1}^n \log_q [f_\epsilon(x_i; \boldsymbol \theta)].
%\end{equation}
%
 The aim is to find $\hat{\boldsymbol \theta}$, for which

\begin{equation}
\hat{\theta} = \arg \max_{\boldsymbol \theta} \left(\psi_\Lambda(f_0,f_1) \right),
\end{equation}
where
\begin{equation}
\psi_\Lambda(f_0,f_1) = \frac{\mathrm{d} \rho_{\Lambda}(f_\epsilon)}{\mathrm{d} \epsilon}\left|_{\epsilon=0} \right.\, .
\end{equation}

The operator $\frac{\mathrm{d}}{\mathrm{d} \epsilon}|_{\epsilon=0}$ describes the change of estimating function $\rho_\Lambda$ under a small contamination of $f_0$ by $f_1$. Note that the derivative can be understood a special case of variational calculus \cite{Gelfand.Fomin.1963} w.r.t parameter $\epsilon$, which can be rewritten as
%\begin{equation}
%\frac{\mathrm{d}}{\mathrm{d} \epsilon}|_{\epsilon=0} \rho_\Lambda(f_\epsilon) = \sum_i \frac{\partial}{\partial f_\epsilon(x_i,\theta)}|_{f_\epsilon = f_0} \frac{\mathrm{d} f_\epsilon(x_i,\theta)}{\mathrm{d} \epsilon}|_{\epsilon=0} = \sum_i \frac{\partial \Lambda(f_\epsilon)}{\partial f_\epsilon(x_i; \boldsymbol \theta)}|_{f_\epsilon = f_0} (f_1(x_i; \boldsymbol \theta) - f_0(x_i; \boldsymbol \theta)) .\end{equation}

\begin{equation}
\psi_\Lambda(f_0,f_1) = \frac{\mathrm{d} \rho_{\Lambda}(f_\epsilon)}{\mathrm{d} \epsilon}\left|_{\epsilon=0} \right. = \sum_{i=1}^n \frac{\partial \Lambda(f_\epsilon(x_i,\boldsymbol \theta)}{\partial f_\epsilon(x_i,\boldsymbol \theta)}|_{f_\epsilon=f_0} \cdot \frac{\ud f_\epsilon(x_i,\boldsymbol \theta)}{\ud \epsilon}|_{\epsilon = 0} = \sum_{i=1}^n \frac{\partial \Lambda(f_\epsilon)}{\partial f_\epsilon(x_i; \boldsymbol \theta)}|_{f_\epsilon = f_0} \cdot (f_1(x_i; \boldsymbol \theta) - f_0(x_i; \boldsymbol \theta)) .\end{equation}

% \sum_i \frac{\partial \Lambda(f_\epsilon(x_i,\boldsymbol \theta)}{\partial f_\epsilon(x_i,\theta)}|_{f_\epsilon=f_0}

% The contamination function $f_1$ is from the eliminated via $\epsilon=0$.

%\Lambda$ represents  generalized logarithm ($\log_q$) from generalized entropy (Tsallis) and $\log$.  $\Lambda$ has to be a strictly monotonic to make a one to one mapping from $f_\epsilon(x;\boldsymbol \theta)$ to $\Lambda$ (see subsection \ref{maxminSMono}).

 %minimizing or

%\noindent $\rho_{*}$ is estimating functions.

 %$f_\epsilon( x; \hat{\boldsymbol \theta})$  converge in probability to $ f_\epsilon(x;\boldsymbol \theta)$ and
 %$f_\epsilon( x; \hat{\boldsymbol \theta})$  is a sampling version of $f_\epsilon(x;\boldsymbol \theta)$ \cite{gutgraduate}. %$\rho_{*}(f_0( x; \hat{\boldsymbol \theta}),f_1(x; \hat{\boldsymbol \theta}))$ is the convergence to $\rho_{*}(f_0(x;\boldsymbol \theta),f_1(x;\boldsymbol \theta))$ and $\rho_{*}(f_0( x; \hat{\boldsymbol \theta}),f_1(x;\hat{\boldsymbol \theta}))$ is a sampling version of $\rho_{*}(f_0(x;\boldsymbol \theta),f_1(x;\boldsymbol \theta))$.
 %It is possible to propose infinite estimating functions from PDFs to fit data. However,  since $\rho_{*}(f_0( x; \hat{\boldsymbol \theta}),f_1(x;\hat{\boldsymbol \theta}))$  comes from only one $\rho_{*}(f_0(x;\boldsymbol \theta),f_1(x;\boldsymbol \theta))$, some PDFs or their convex combination from LIDs in $\Lambda$ are capable of fitting data. Expecting that all of them gives the best fitting is not true.

%\subsection{Property of Function $ \Lambda $}\label{maxminSMono}

Generally, $\log_q$ substituted to $\Lambda$ is not the only possibility for generalization of MLE. The function $\Lambda$ in equation \eqref{estimatingPDF} has to be a strictly monotonic function of its argument. This can be easily investigated by the first derivative of $\Lambda$. Let us focus on $\Lambda(f_\epsilon)=\log_q(f_\epsilon)$. There are no roots of first derivatives for $\log_q(f_\epsilon)$.  $\Lambda$ is an one to one mapping, but there is no a zero value coming from $\Lambda$ except $\log(1)=0$ and $\log_q(1)=0$. The second derivative is $-qf_\epsilon^{-q-1}$. For $q>0$, the second derivative test shows that these functions are concave. Since $\log_q$ is concave for $q>0$ and also there are no roots for the first derivative, except to $f_\epsilon=1$ which is impossible, they can satisfy to be strictly monotonic functions. As a result, this is an important property to use this function for estimations of parameters in a PDF.

\subsection{Estimating functions from $f_0$ and Least informative distributions in $\log_q$ functions}
%Generalized exponential and $\log$ functions are introduced by \cite{Kaniadakisetal04,Kaniadakisetal05,MatsuWadaGex,kaniadakissc08,Mittal75}.
%An  estimating function from $\log_{q}$ is  \cite{tsallis1988possible}:
Let us now focus on LIDs obtained from $q$-log-likelihood function $\rho_{\log_q}$ defined as
\begin{equation}\label{eqlnq}
 \rho_{\log_{q}}(f_0)  =  \sum_{i=1}^n \frac{f_0(x_i;\boldsymbol \theta)^{1-q}-1}{1-q}.
\end{equation}
Now, let us consider a convex combination $f_\epsilon$ and plug it into $\rho_{\log_q}$: % is used to produce ,
\begin{equation}\label{LIDTsallisstart}
 \rho_{\log_q}(f_\epsilon) = \frac{1}{1-q} \bigg[\sum_{i=1}^{n}  [(1-\epsilon)f_0(x_i;\boldsymbol \theta)  + \epsilon f_1(x_i;\boldsymbol \theta)]^{1-q} -1  \bigg],
 \end{equation}
after taking derivative  w.r.t $\epsilon$ and putting $\epsilon=0$ , we get % setting $\epsilon$ to be zero gives

\begin{equation}\label{eqlnqLID}
 \psi_{\log_{q}}(f_0,f_1) =   \sum_{i=1}^n f_0(x_i;\boldsymbol \theta)^{-q} [f_1(x_i;\boldsymbol \theta) - f_0(x_i;\boldsymbol \theta)],
\end{equation}
\noindent where $q$  is  a tuning parameter for robust estimation to produce functions that are neighborhood for $f_{\epsilon}$. The equation \eqref{eqlnqLID} is defined to be a LID. It can be rewritten as a similar form to mixed distribution

\begin{equation}\label{eqlnqLIDWeight}
 \psi_{\log_{q}}(f_0,f_1) =   \sum_{i=1}^n [ w_1(x_i;\boldsymbol \theta) f_1(x_i;\boldsymbol \theta) + w_0(x_i;\boldsymbol \theta) f_0(x_i;\boldsymbol \theta)],
\end{equation}
\noindent where $w_1=f_0(x_i;\boldsymbol \theta)^{-q}$ and $w_0=-f_0(x_i;\boldsymbol \theta)^{-q}$.

%\begin{equation}\label{eqlnqLIDWeight}
% \rho_{\log_{q}}(f_0,f_1) =   \sum_{i=1}^n f_0(x_i;\boldsymbol \theta)^{-q} f_1(x_i;\boldsymbol \theta) - f_0(x_i;\boldsymbol \theta)^{-q} f_0(x_i;\boldsymbol \theta),
%\end{equation}

As a special case of these generalized functions, LID obtained from MLE (i.e., when $\Lambda = \log$) is
\begin{equation}\label{eqlnLID}
 \psi_{\log}(f_0,f_1) =   \sum_{i=1}^n f_0(x_i;\boldsymbol \theta)^{-1} [f_1(x_i;\boldsymbol \theta) - f_0(x_i;\boldsymbol \theta)].
\end{equation}
\noindent  LID from equation  \eqref{eqlnLID} can also be rewritten as $\sum_{i=1}^n [\frac{f_1(x_i;\boldsymbol \theta)}{f_0(x_i;\boldsymbol \theta)} - 1]$.  %A comment should be given to show role of $(1-\epsilon)f_0  + \epsilon f_1$ used in LID. Let us consider that a value of $f_0$ is bigger than a value of $f_1$, then a value of $f_1/f_0-1$ as high will support that an observation will be a member of $f_0$. When a value of $f_1$ is bigger than a value of $f_0$, a value of $f_1/f_0-1$ is low, which is a good result for us to imply that an observation comes from $f_1$. In this case, a low value of $f_0$ is supported by a low value of  $f_1/f_0-1$ as well.
 Having $f_0$ and $f_1$ together means that data are distributed non identically. It can also be considered that equations \eqref{eqlnqLID} and \eqref{eqlnLID} are mixing distributions, i.e, $f_\epsilon$.

One can consider that generalized entropies, generalized logarithms (generalized exponentials) and divergences \cite{JizbaKorbel16,MachadoEntropy,UbriacoFC,Wadatwopara,HanThur11} and references therein can be applied to get new $\psi$, but these functions and function $\Lambda$ play same role as a weighted MLE form for estimations of parameters. In this context, we can have functions that can be equal to each other for a value of parameter. For example, $\frac{1}{1-q}(f_0^{1-q}-1)$ and $\frac{1}{1-q} (f_0^q-1)$ are equal to each other when $q=1/2$.

\section{Information criterions based on estimating functions $\rho_{\log}$ from MLE,  $\rho_{\log_q}$ from MqLE, $\psi$ from Huber M-Estimation and $\psi_{\log}$, $\psi_{\log_q}$ from LID}\label{modAIC}
Information criterion (IC) is a tool for assessing the fitting performance of functions. Different tools are proposed by \cite{Aka73,Hamparsum87}. After proposing estimating functions from $\Lambda$, we will have another problem to test the fitting performance of $\rho_{\log}$ from MLE,  $\rho_{\log_q}$ from MqLE, $\psi$ from Huber M-Estimation and $\psi_{\log}$, $\psi_{\log_q}$ from LID. For this aim, robust information criterion (RIC) formulae are used to determine value of tuning parameter $q$. % in operator $D^{q}$.

Let us consider the equation \eqref{estimatingPDF} including the function $\log$ as a special case of $\Lambda$. Due to this reason, IC is
\begin{equation}\label{AICformula}
IC(f_0,c_k)=-2\rho_{\log}(f_0) + c_k.
\end{equation}
IC can have two forms that are AIC and BIC. Performance of AIC depends on  penalty term $c_k=2k$, which is a deficiency of AIC  \cite{Hamparsum87} and references therein. As an alternative to AIC, BIC was proposed when $c_k=\log(n)k$. We propose robust version of ICs via replacing  $\rho_{\log}$ with $\psi_{\log_q}$. Robust versions of ICs have been proposed by \cite{FerrariOpere,Zhang10,RoncSS97} as a same approach we proposed here.  ICs can be considered as an appropriate form for  $\psi_{\log_q}$ from LID
\begin{equation}\label{MAICalphabeta}
  RIC_q(f_0,f_1,c_k)= -2 \psi_{\log_q}(f_0,f_1)  + c_k,
\end{equation}
\noindent where $k$ is number of estimated parameters.  When $c_k$ is $2k$ and $\log(n)k$, ICs are robust Akaike and robust Bayesian ICs ($RAIC_{q}$ and $RBIC_{q}$), respectively. $RIC_{q}$ for $f_0$ and $\psi$ from Huber M-Estimation can be rewritten in the similar way. MLE is maximization of $\rho_{\log}(f_0)$ according to parameters. If this function has a maximum value, then $-2\rho_{\log}(f_0)$ and $-2 \psi_{\log_q}$ will have a minimum value. Minimum values of equations \eqref{AICformula}-\eqref{MAICalphabeta} mean that fitting performance is accomplished \cite{Aka73,Hamparsum87,FerrariOpere,Zhang10,RoncSS97}.

\vspace{3mm}

\section{Real data analyzing procedure and artificial data generated from underlying distribution $f_0$}\label{realsection} %
Optimizing  $\psi_{\log_q}$ defined by  $\Lambda$ according to parameters in PDFs $f_0$ and $f_1$ produces M-estimators $\hat{a}$ and $\hat{b}$ from LID
\begin{equation}\label{Malpbetab}
(\hat{a}_{\psi},\hat{b}_{\psi}) := \underset{a ~ \text{and} ~ b}{\arg  \max}  ~~~ \psi_{\log_q}(f_0,f_1).
\end{equation}
If only $f_0$ is chosen for $\rho_{\log_q}$, then estimators $\hat{a}$ and $\hat{b}$ will be obtained from a PDF. Since $\psi_{\log_q}$ and $\rho_{\log_q}$ are nonlinear functions according to the parameters in a PDF, an optimization method is essential to use. The hybrid genetic algorithm (HGA) in MATLAB R2016a was used to get the estimates of parameters. Due to working principle of HGA, the prescribed interval for parameters $a$ and $b$ is $[0,50]$.

The possible smallest value of $RIC_{q}$ and the frequency (counted data at divided intervals of domain) are accepted to be a best choice among them while performing our trying for different values of parameter $q$, $p_0$ in $f_0$ and $p_1$ in $f_1$ as shape parameters. Since simultaneous estimations of the shape parameters $p_0$ and $p_1$, the parameter $q$ in $\log_q$ and other parameters, such as $a$, $b$, $\mu$ and $\sigma$ in PDFs are not easy and also the parameters $p_0$, $p_1$ and $q$ are known to be the tuning parameters for robust estimations of parameters in PDFs \cite{Hub64,ArslanGenc09} and references therein, we consult to use the information criterions $RIC_{q}$ and the frequency in order to adjust the values of tuning parameters. For these reasons, the parameters $p_0$, $p_1$ and $q$ are taken to be fixed while performing the estimation procedure of the parameters $a$, $b$ with fixed $q$ from distributions  on $[0,\infty)$: see the examples in subsections \ref{ex1abPDFs} and \ref{ex2abPDFs} for estimations of shape $a$, scale $b$ parameters and  $\mu$, $\sigma$ with fixed $q$, $p_0$ and $p_1$ from distributions on  $(-\infty,\infty)$: see the examples in \ref{appBlocsca} for estimations of location $\mu$ and scale $\sigma$ parameters.

Each estimating function is comparable in itself, because they are different functions from each others. However, $\rho_{\log}(f_0)$  is comparable, because they set same mapping from PDF to $\log[f_0(x;\boldsymbol \theta)]$. All of these estimating functions are different from each others, because they have not only different analytical form but also different values for $q$. In this context, it is very difficult to know the best function for fitting, therefore we need to propose functions that are neighborhood to each others, which can help us to perform the best fitting on non identically distributed data as well.

An outlier that leads to have the non identical case is obtained by maximum value $x$ multiplied by 2, that is $2\text{max}(x)$, for both of two examples. We added one outlier in the left and right side of data set for examples in \ref{ex1musig} and \ref{ex2musig}. Thus, we can test efficiency and robustness of estimators obtained by $\rho_{\log}$ from MLE,  $\rho_{\log_q}$ from MqLE, $\psi$ from Huber M-Estimation and $\psi_{\log}$, $\psi_{\log_q}$ from LID.

%Note that we focus on the non-outlier case, because we keep the property of data set. In other words, we choose the non-distorted data set. We follow this procedure for other examples throughout the text.
\vspace{3mm}
\subsection{Example 1: M-Estimations and MLEs of the parameters $a$ and $b$ in PDFs}\label{ex1abPDFs}

Observations from temperature of 1951-2017 years \cite{Datasetref} are used in analyzing procedure, because there are important changes in temperature on Earth. We are interested in the estimations of parameters $a$ and $b$ in PDFs at \ref{f0collection}. The value of outlier is $6.5973$.

\begin{table}[!htb]
\centering
\caption{Estimates of parameters $a$ and $b$ from different estimating functions without and with one outlier for temperature data}
\label{ex1ab2}
\scalebox{0.99}{
\begin{tabular}{ccccc}
Estimating Functions & $\hat{a}$  & $\hat{b}$ & $RAIC_{q}$   & $RBIC_{q}$    \\ \hline \hline
$\psi_{\log_{q=.007}}(f_0=\text{Gamma},f_1=\text{Weibull})$  &       2.9699         &   0.2431     &      56.6114    &   61.6110   \\
One Outlier &      2.9699 &  0.2431   &     56.6114    &  61.6331   \\	 \hline
$\rho_{\log_{q=0.53}}(f_0=\text{Gamma})$  &       3.0624      &  0.2422     &110.9574     &115.9570  \\
One Outlier &    3.0402      &  0.2395    &     115.3203    &  120.3420   \\	 \hline \hline
Estimating Function & $\hat{a}$  & $\hat{b}$ & AIC   & BIC    \\ \hline \hline
$\rho_{\log}$(Gamma)  & 1.9280 & 0.5920 & 188.8185 &193.8182  \\
One Outlier &1.7610  &0.6822  &204.6059  & 209.6276 \\	\hline
$\rho_{\log}$(Weibull) & 1.3874 & 1.2603 & 193.2213  & 198.2210  \\
One Outlier & 1.2882 & 1.3097 & 209.8081 & 214.8298 \\	\hline
$\rho_{\log}$(Burr) & 2.4977 & 0.8109 & 185.1382 & 190.1378 \\
One Outlier &2.3918  & 0.8416 & 196.7629 & 201.7847 \\	\hline\hline
\end{tabular}}
\end{table}
\begin{figure}[!htb]
\centering
 \begin{subfigure}{.67\linewidth}
    \includegraphics[width=1.05\textwidth]{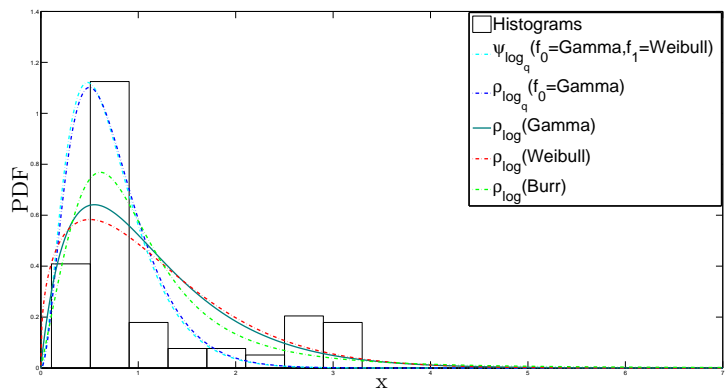}
    \caption{M-estimators and MLE of parameters in $f_0$ and LID}
  \end{subfigure}
 \begin{subfigure}{.67\linewidth}
    \includegraphics[width=1.05\textwidth]{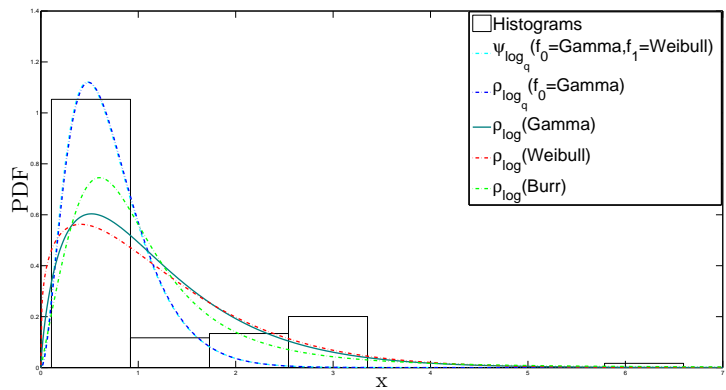}  %
 \caption{M-estimators and MLE of parameters in $f_0$ and LID when one outlier is added} %
  \end{subfigure}
  \caption{M-estimators and MLE of parameters in $f_0$ and LID for temperature data}
  \label{Ex1figtoplucaAB}
  \end{figure}

   The  values of $RAIC_{q}$ and $RBIC_{q}$ for $\psi_{\log_{q=.007}}(f_0=\text{Gamma},f_1=\text{Weibull})$ and $\rho_{\log_{q=.53}}(f_0=\text{Gamma})$ are not comparable, because the mapping region is not same. For this reason, $RAIC_{q}$ and   $RBIC_{q}$ for these estimating functions are not homogeneous to compare in their self. Therefore, we need to take in account the counted data at divided intervals of domain.   The estimates from $\psi_{\log_{q=.007}}(f_0=\text{Gamma},f_1=\text{Weibull})$ show more robustness when they are compared with  $\rho_{\log_{q=.53}}(f_0=\text{Gamma})$.

   Among $\rho_{\log}$ functions, $\rho_{\log}$(Burr) has the smallest AIC and BIC values and $\rho_{\log}$(Gamma) is second function on fitting performance for non-outlier case. The outlier case of Burr III distribution from $\rho_{\log}$ has the smallest values of AIC and BIC among information criterions. Since Burr III distribution is heavy tailed, it can generate data that are far from the underlying distribution. Therefore, Gamma distribution was preferred as an underlying distribution.  We can also get the efficient fitting when Gamma is underlying and Weibull is contamination distributions.

\vspace{2mm}
\subsection{Example 2: M-Estimations and MLEs of the parameters $a$ and $b$ in PDFs}\label{ex2abPDFs}

One can get data from page https://legacy.bas.ac.uk/met/READER/surface/Grytviken.All.temperature.txt:Grytviken temperature in November at each year from 1905 and 2017. The data in some years are missed. The value of outlier is $13.2$.

\begin{table}[!htb]
\centering
\caption{Estimates of parameters $a$ and $b$ from different estimating functions without and with one outlier for temperature data}
\label{ex1ab2}
\scalebox{0.99}{
\begin{tabular}{ccccc}
Estimating Functions & $\hat{a}$  & $\hat{b}$ & $RAIC_{q}$   & $RBIC_{q}$    \\ \hline \hline
$\psi_{\log_{q=.61}}(f_0=\text{Gamma},f_1=\text{Weibull})$  &   7.5454  &   0.4389   &   106.8688 & 111.9765     \\
One Outlier &   7.5234  &  0.4491  &     105.4343 & 110.5630 \\	 \hline
$\rho_{\log_{q=.22}}(f_0=\text{Gamma})$  &        5.4602 &   0.5590      & 173.6961 & 178.8039   \\
One Outlier &   5.7471 & 0.5201    &     175.5467 &  180.6754  \\	 \hline \hline
Estimating Function & $\hat{a}$  & $\hat{b}$ & AIC   & BIC    \\ \hline \hline
$\rho_{\log}$(Gamma)  &     4.1357  &  0.7868  &     325.1162 & 330.2239  \\
One Outlier &     4.0627  &  0.8122  &   344.6509  & 349.7796    \\	\hline
$\rho_{\log}$(Weibull) &       2.7380  &  3.5414    &   313.6991 &  318.8069  \\
One Outlier &     2.0953 &   3.6720 &   353.0678 &  358.1965 \\	\hline
$\rho_{\log}$(Burr) &     1.9911 &   5.5948 &   366.7043 &  371.8121 \\
One Outlier &    1.9664 &   5.5872 &   377.3481 &  382.4768 \\	\hline\hline
\end{tabular}}
\end{table}
\begin{figure}[!htb]
\centering
 \begin{subfigure}{.67\linewidth}
    \includegraphics[width=1.05\textwidth]{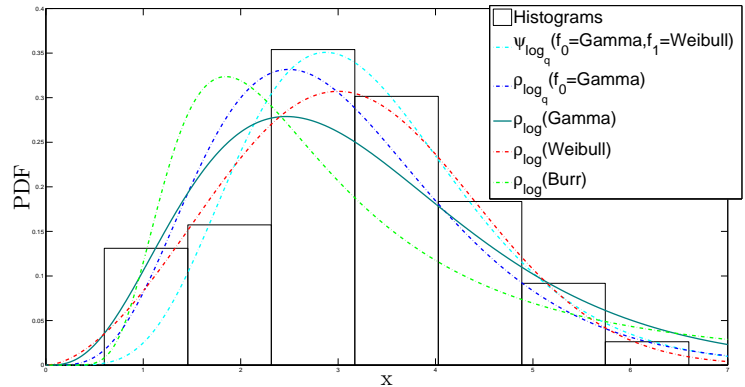}
    \caption{M-estimators and MLE of parameters in $f_0$ and LID}
  \end{subfigure}
 \begin{subfigure}{.67\linewidth}
    \includegraphics[width=1.05\textwidth]{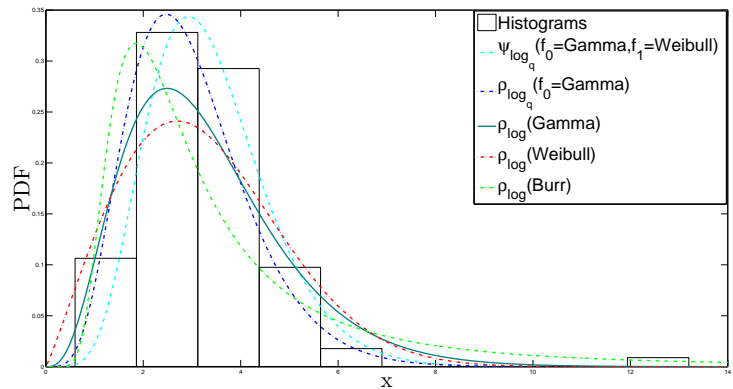}  %
 \caption{M-estimators and MLE of parameters in $f_0$ and LID when one outlier is added} %
  \end{subfigure}
  \caption{M-estimators and MLE of parameters in $f_0$ and LID for temperature data}
  \label{Ex2figtoplucaAB}
  \end{figure}

  The  values of $RAIC_{q}$ and $RBIC_{q}$ for $\psi_{\log_{q=.61}}(f_0=\text{Gamma},f_1=\text{Weibull})$ and $\rho_{\log_{q=.22}}(f_0=\text{Gamma})$ are not comparable, as implied by example 1. The estimates from $\rho_{\log_{q=.61}}(f_0=\text{Gamma},f_1=\text{Weibull})$ show more robustness when they are compared with  $\rho_{\log_{q=.22}}(f_0=\text{Gamma})$. Among $\rho_{\log}$ functions, $\rho_{\log}$(Weibull) has the smallest AIC and BIC values and $\rho_{\log}$(Gamma) is second function on fitting performance for non-outlier case. The outlier case of Burr III distribution from  $\rho_{\log}$ has the highest values of AIC and BIC among information criterions. Since Burr III distribution is heavy tailed and modelling data in example 2 is not good and also it can generate data that are far from the underlying distribution. Therefore, Gamma distribution was preferred as an underlying distribution. We can also get the efficient fitting when Gamma is underlying and Weibull is contamination distributions. As it is seen for both of two examples, the best fitting can be accomplished by such mixed distribution.

For both of examples, the parameter $q$ is responsible for determining the tails in the left and right sides of PDF together with the overall shape of PDF, which means that we can overcome non-identicality problem in a data. Thus, the robustness and efficiency can be guaranteed. Data have their nature, therefore the modelling of them depends on the properly chosen values for the parameter $q$ and a function considered for fitting data. In this sense, the information criterions have to be used. When we look at results, $RBIC_{q}$ is sensitive, thus it is more informative when it is compared with $RAIC_{q}$, because $RBIC_{q}$ can detect outlier case especially for example 1. Let us consider about decreasing values of $RAIC_{q}$ and $RBIC_{q}$ in outlier case at example 2. Then, it is observed that the best fitting has been performed in the outlier case, because the counted data of  $\psi_{\log_{q}}(f_0=\text{Gamma},f_1=\text{Weibull})$ for outlier case in Table 4 support that the best fitting has been done by $\psi_{\log_q}(f_0=\text{Gamma},f_1=\text{Weibull})$ very well for the general part of data set. In the decision procedure about fitting performance, we need to draw PDFs of underlying distributions at Figures \ref{Ex1figtoplucaAB} and \ref{Ex2figtoplucaAB} for the estimated values of parameters from $\rho_{\log_q}$, generate artificial data (see the subsection \ref{ArtdatasetEx1} for details) and have information criterions as well.

Let us give a comment for Burr III distribution in $\rho_{\log}$ from two examples \ref{ex1abPDFs} and \ref{ex2abPDFs}. Although Burr III in example  \ref{ex1abPDFs} has a good fitting among $\rho_{\log}$ functions, i.e. MLE,  Burr III in example  \ref{ex2abPDFs} does not have a good fitting among MLEs. As it is seen, the data and $\rho_{\log}(f_0)$ have to accommodate  each others if we get the more precise estimators for the parameters. In LID and MqLE case, it is possible to drive the function via parameter $q$ in $\log_q$. Thus, we can make an accommodation between data and function.

\vspace{5mm}

\subsection{Simulation: Artificial data set generated from distribution  $f_0$ on $ \mathbb{S} = [0,\infty)$}\label{ArtdatasetEx1}

For case of estimating the parameters $a$ and $b$ in a PDF at \ref{f0collection}, the underlying $f_0$ is used to generate artificial data from Gamma and Weibull distribution and draw PDFs with the estimated values from $\rho_{\log_q}$, resp. $\psi_{\log_q}$.  The histograms in Figures \ref{Ex1figtoplucaAB} and \ref{Ex2figtoplucaAB} are given to illustrate the similarity to  histograms of real data.  Since we are interested in the estimations of shape and scale parameters, the artificial data must be generated to see the behaviour of shape generally.  The number of replication is 100 000. Data set generated from underlying distribution is sorted in each $n$. After sorting data in each $n$, an arithmetic mean of 100 000 artificial data is obtained for the sample size $n=90$ of example 1 in the case of estimating the parameters $a$ and $b$.  After performing simulation, we can have more precise decision about which a PDF with its estimated values of parameters from $\rho_{\log}$, $\rho_{\log_q}$ and  $\psi_{\log_q}$ is appropriate to represent real data very precisely. Each of function in Table \ref{Ex1Tablecount} has different counted data at bin ranges $[0,0.5,1.5,2.5,20]$. $\rho_{\log}$(Gamma), $\rho_{\log}$(Weibull) and $\rho_{\log}$(Burr) tend to fit the data in tail. However,     $\psi_{\log_q}(f_0=\text{Gamma},f_1=\text{Weibull})$ and $\rho_{\log_q}(f_0=\text{Gamma})$ are resistant to data in tail. Additionally, the numbers 52+7=59 from real data and 54+5=59 of artificial data from $\psi_{\log_q}(f_0=\text{Gamma},f_1=\text{Weibull})$ are same to each others  when the counted data at left side are ignored.

\begin{table}[!htb]
\centering
\caption{The counted data at each bins $[0,0.5,1.5,2.5,20]$}
\label{Simabex1}
\scalebox{1.15}{
\label{Ex1Tablecount}
\begin{tabular}{c|ccccc} \hline \hline
Real data &    16  &  52   &  7  &  15 &    0  \\
One Outlier &  16   & 52  &   7  &  16  &   0 \\  \hline
$\psi_{\log_q}(f_0=\text{Gamma},f_1=\text{Weibull})$ &  31  &   54 &    5 &    0 &    0 \\
One Outlier &    31 &   55   &  5 &    0  &   0  \\ \hline
$\rho_{\log_q}(f_0=\text{Gamma})$ &  29  &  56   &  5  &   0 &    0 \\
One Outlier &  31  &  55   &  5   &  0   &  0  \\ \hline
$\rho_{\log}$(Gamma)  &   20 &   46  &  18 &    6   &  0 \\
One Outlier &   20   & 45 &   18  &   8  &   0  \\ \hline
$\rho_{\log}$(Weibull) &  22 &   43  &  18  &   7  &   0  \\
One Outlier & 23   & 40  &  19   &  9  &   0 \\ \hline
$\rho_{\log}$(Burr) &   19  &  51  &  13 &    7 &    0  \\
One Outlier &   19 &   50 &   14  &   8  &   0 \\  \hline
\end{tabular}}
\end{table}

\begin{table}[!htb]
\centering
\caption{The counted data at each bins $[0,1,2,3,4,5,50]$}
\scalebox{1.12}{
\label{Ex2Tablecount}
\begin{tabular}{c|ccccccc} \hline \hline
Real data &         6  &  11  &  23 &   30  &  16   &  9   &    0  \\
One Outlier &       6  &  11  & 23  &  30   &  16   &  10  &   0 \\  \hline
$\psi_{\log_q}(f_0=\text{Gamma},f_1=\text{Weibull})$ &       0  &  12 &   30 &   29 &   15 &    9 &    0  \\
One Outlier &        0   & 11 &   29  &  30  &  16  &  10   &  0   \\ \hline
$\rho_{\log_q}(f_0=\text{Gamma})$ &      2 &   19   & 30    &24  &  12&     8  &   0   \\
One Outlier &       2   & 19 &   33  &  24 &   11  &   7   &  0   \\ \hline
$\rho_{\log}$(Gamma)  &      3    &19   & 26 &   21  &  13 &   13&     0   \\
One Outlier &        3   & 18  &  26  &  21  &  14 &   14   &  0     \\ \hline
$\rho_{\log}$(Weibull) &  3   & 15  &  27   & 27   & 16  &   7   &  0   \\
One Outlier &      6 &   17&    23 &   21  &  15&    14 &    0    \\ \hline
$\rho_{\log}$(Burr) &       2 &   25  &  25  &  15  &   9   & 19   &  0   \\
One Outlier &      2  &  25 &   25    &15&     9  &  20 &    0   \\  \hline
\end{tabular}}
\end{table}

Each of function in Table \ref{Ex2Tablecount} has different counted data at bin ranges $[0,0.5,1.5,2.5,50]$. Generating the artificial data is performed for the sample size $n=95$ at  100 000 replications for example 2 in the case of estimating the parameters $a$ and $b$. $\rho_{\log}$(Gamma), $\rho_{\log}$(Weibull) and $\rho_{\log}$(Burr) tend to fit the data in tail especially in the case of an outlier. However, $\psi_{\log_q}(f_0=\text{Gamma},f_1=\text{Weibull})$ is resistant to data in left tail and fits the right side of data well (see also Figure \ref{Ex2figtoplucaAB}-(a) and (b)). Additionally, it can be observed that the general part of data can be represented by LID case. Especially, it can represent real data in outlier case very well when it is compared with functions $\rho_{\log_q}(f_0=\text{Gamma})$, $\rho_{\log}$(Gamma), $\rho_{\log}$(Weibull) and $\rho_{\log}$(Burr).

\vspace{15mm}

\section{Conclusions}\label{conclusionsection}

Estimating functions from LIDs have been proposed when  the data are composed of two PDFs. Since we use convex combination of two functions, the more informative data analyzing procedure can be done. We used $\Lambda$ to propose new estimating functions $\psi_{\log_q}$ and $\psi_{\log}$. Although we eliminate $f_1$ from $f_\epsilon$, it is interesting that the role of $f_1$ keeps in LID. Thus, this procedure is also considered as a mixing of two PDFs $f_0$ and $f_1$. The contamination rate $\epsilon$ is not known exactly or getting the estimations of mixing proportions is difficult task. Using $\psi_{\log_q}$ is better than using $\psi_{\log}$ to estimate parameters efficiently and robustly, because the parameter $q$ and also LID  in $\log_q$ are advantages for us to propose a flexible function. This flexibility produces efficient and robust estimators from the neighborhood of $f_0(x;\boldsymbol \theta)$ and especially $f_{\epsilon}(x;\boldsymbol \theta)$  in $\log_q$.  $RIC_{q}$ is proposed to assess the fitting performance of estimating functions from $\log_q$. Thus, the value of tuning parameter $q$ in $\log_q$ can be determined according to minimum values of $RIC_{q}$. Here, generating the artificial data is also an another important issue in order to determine the value of tuning parameter $q$. Real data sets were provided to show the fitting competence of estimating functions. In the illustrating the performance of M-estimation, it is observed that $\psi_{\log_q}(f_0,f_1)$ case can have values for $q$ at a small range when it is compared with that of $\rho_{\log_q}(f_0)$. In this sense, $\psi_{\log_q}(f_0,f_1)$ can be an advantage to determine the value of $q$ easily, as it is observed from  example 1 for the estimations of the parameters $a$ and $b$ and also both of two examples in the estimations of the parameters $\mu$ and $\sigma$. %as it is observed for a such kind of data sets,

As it is well known, the concavity  property of $\Lambda$ is important to make an one to one mapping from $f_{\epsilon}(x;\boldsymbol \theta)$ to $\Lambda$. Otherwise, $\Lambda$ does not give correct result for mapping.  One can use the proposed $\psi_{\log_q}$ to produce  LIDs by choosing arbitrary $f_0$ and $f_1$ if parameters in PDFs are same property, such as shape, scale and location. Thus, we have flexible LID that is a convex combination of $f_0$ and $f_1$ in $\Lambda$.  In future, we will prepare a package for $\psi_{\log_q}$ in univariate and multivariate variables at open access R software. The regression case can be done as an application of location model. The robust test statistics based on $\log_q$ are our ongoing research \cite{CanKor16,Yalcinetal}. We will also study the score functions for these estimating functions and its connection with Fisher metric \cite{CanKor16}. Many phenomena (data in signal \cite{MatsuHen} and image \cite{Peters} processings, climate change, medical issues, etc.) which are modelled by the parametric models can be analysed by means of this package. Thus, precise and robust estimations of parameters can be done via $\log_q$ and LID based on $\log_q$.
%Trying other generalized logarithms for the sake of determining a small interval for changing of tuning parameter in generalized logarithms if it so is an open problem.
\section{Acknowledgements}
%%%%%%%%%%%%%%%%%%%%%%%%%%%%%%%%%%%%%%%%%%%%%%%%%%%%%%%%%%%%%%%%%%%%%%%%
We are indebted to Prof. Dr. James F. Peters from the University of Manitoba, Winnipeg, Canada, for critical reading. We also thank to Foreign Language School of  U\c{s}ak University for editing language and partial support  from Turkish government for M.N.\c{C}. J. K. was supported by the Austrian Science Fund, grant No. I 3073-N32 and by the Czech Science Foundation, grant No. 17-33812L.
\appendix

%\newpage

\section{PDFs used to get estimating functions $\rho_{\log}$ from MLE,  $\rho_{\log_q}$ from MqLE, $\psi$ from Huber M-Estimation and $\psi_{\log}$, $\psi_{\log_q}$ from LID}\label{f0collection}
%The distributions are considered to get estimating functions $\rho_{*}$ for only $f_0$ and LID.
\begin{table}[!htb]
\begin{center}%Exponential
\caption{Distributions on  positive half-plane \cite{Weibuliseverywhere,Burr42} (Weibull, Gamma and Burr type III) and on whole real axis \cite{ArslanGenc09,Cankaya} (Exponential power and Generalized t). PDFs include gamma $\Gamma(\cdot)$ and beta $B(\cdot,\cdot)$ functions. }
\scalebox{0.93}{
\label{PDFstable}
\begin{tabular}{ccc}
\hline
%\multicolumn{3}{c}{Distributions}
Distributions & $ \boldsymbol \theta $: Role of Parameters & PDF  \\ \hline
Half-plane: $ \mathbb{S} = [0,\infty)$ &  &   \\ \hline
%%%%%%%%%%%%%%%%%%%%%%%%%%%%%%%%%
Weibull &
\begin{tabular}{c}  $a$: Shape,
$b$: Scale \end{tabular}
& $f_{\text{Weibull}}(x;a,b)=\frac{a}{b} ( \frac{x}{b} )^{a-1} exp\{- (\frac{x}{b} )^a \}$\\
%%%%%%%%%%%%%%%%%%%%%%%%%%%%%%%%%
\begin{tabular}{c}Gamma \end{tabular}
&\begin{tabular}{c}  $a$: Shape,
$b$: Scale \end{tabular}
& $f_{\text{Gamma}}(x;a,b)=\frac{1}{\Gamma(a)b^a}x^{a-1}exp\{-\frac{x}{b} \}$ \\
%%%%%%%%%%%%%%%%%%%%%%%%%%%%%%%%%%
%Exponential&  $a$: Scale &  $f_{\text{Exp}}(x;a)=\frac{1}{a} exp\{-\frac{x}{a}\}$ \\
%%%%%%%%%%%%%%%%%%%%%%%%%%%%%%%%%%
Burr type III& \begin{tabular}{c} $a,b$: Shape
 \end{tabular}
&  $f_{\text{Burr}}(x;a,b)=abx^{-(a+1)}(1+x^{-a})^{-(b+1)}$  \\
%%%%%%%%%%%%%%%%%%%%%%%%%%%%%%%%%%
\hline
Whole plane: $ \mathbb{S} = (-\infty,\infty)$ &  &  \\ \hline
\hline
%%%%%%%%%%%%%%%%%%%%%%%%%%%%%%%%%
Exponential power & \begin{tabular}{c}
$p,\sigma,\eta>0$\\
$p$: Shape for peakedness,
$\sigma$: Scale\\
$\eta$: Nuisance,
$\mu \in \mathbb{R}$: Location
\end{tabular}  & $f_{\text{EP}}(x;\mu,\sigma,p,\eta)=\frac{p}{2\sigma\eta^{1/p}\Gamma(\frac{1}{p})}exp\{-(\frac{|x - \mu|}{\eta{^{1/p}}\sigma} )^{p} \}$ \\
\hline
%%%%%%%%%%%%%%%%%%%%%%%%%%%%%%%%%%
Generalized t &
\begin{tabular}{c}
$p,\sigma,\nu>0$ \\
$p$: Shape for peakedness,
$\sigma$: Scale   \\
$\nu$: Shape for tail,
$\mu \in \mathbb{R}$: Location
\end{tabular}
& $f_{\text{Gt}}(x;\mu,\sigma,p,\nu)=\frac{p}{2B(1/p,\nu)\nu^{1/p}\sigma}(1+(\frac{|x-\mu|}{\nu^{1/p}\sigma})^p)^{-(\nu+1/p)}$ \\
\hline
\end{tabular}}
\end{center}
\end{table}

\section{Two Examples for M-estimations and MLEs of location $\mu$ and scale $\sigma$ parameters}\label{appBlocsca}
%\newpage

\subsection{Example 1}\label{ex1musig}

A data set is  NCI60 cancer cell line panel. A protein data coded as ME:UACC-257 from Lysate Array at a website https://discover.nci.nih.gov/cellminer/ is analysed. The parameters $\mu$ and $\sigma$ are estimated by using different estimating functions to see tendency (location) and spread (scale) of protein in cancer cell. The maximum value as an outlier is 11.642 at positive and -11.642 at negative sides of real axis. Therefore, we keep the symmetry of data. In HGA, our prescribed intervals for $\mu$ and $\sigma$ are $[-50,50]$ and $[0,50]$, respectively.

%from equation \eqref{Malpbetab}

We will give the comments of results and some values of $p$ for illustrative purpose on modelling capability, comparison with Huber M-estimation and MLE. Among tried values of parameters $p$ and $q$, $RIC_{q}$ values of them in Table \ref{ex1tablemusigma} are given to see the fitting performance of PDFs of which parameters are estimated by  $\rho_{\log}$ from MLE,  $\rho_{\log_q}$ from MqLE, $\psi$ from Huber M-Estimation and $\psi_{\log}$, $\psi_{\log_q}$ from LID. 

\begin{figure}[!]
\centering
   \begin{subfigure}{.67\linewidth}
    \includegraphics[width=0.82\textwidth]{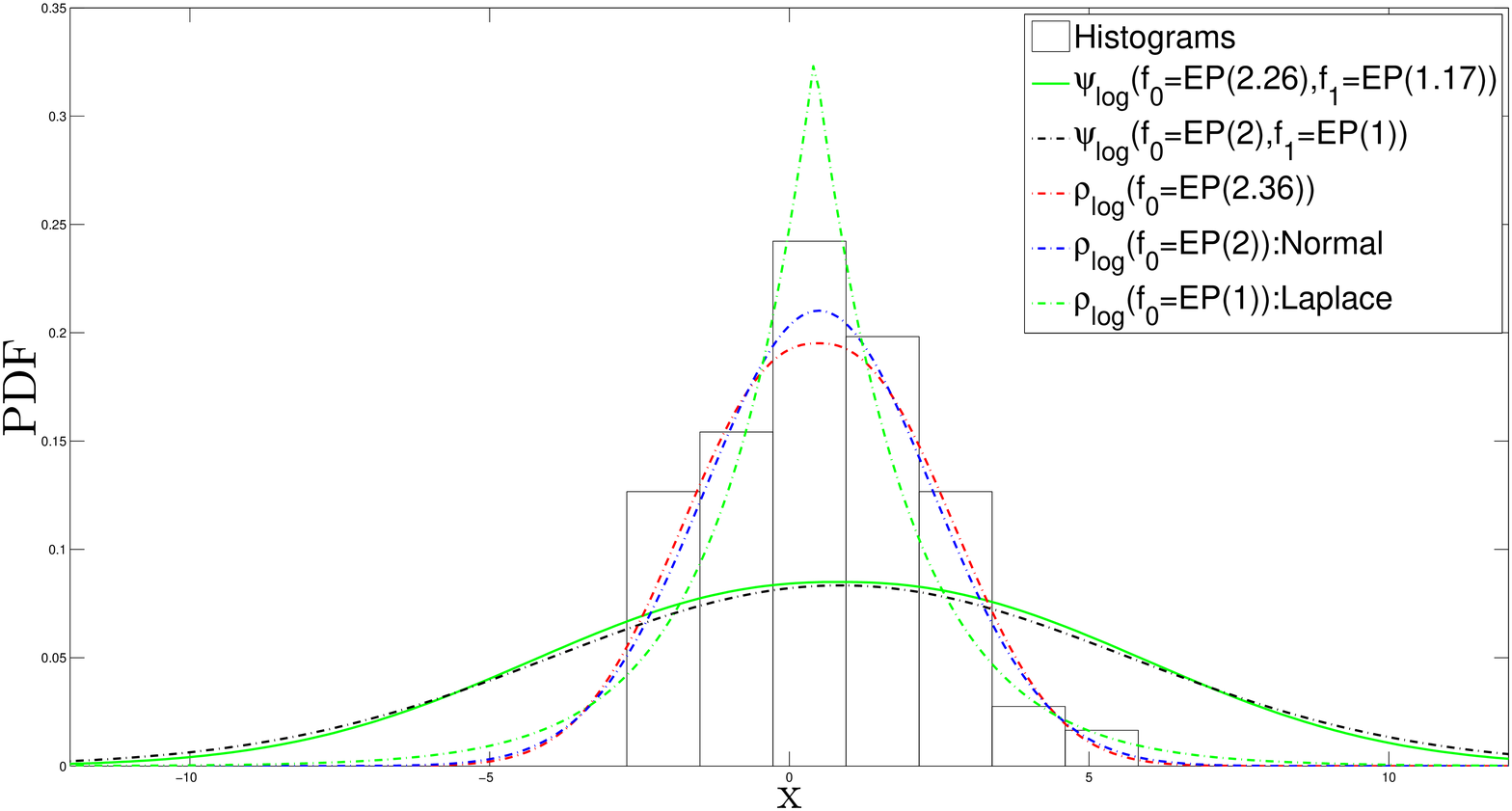}%pdfa2a1collectionDigerleri
    \caption{Non robust M-estimators and MLE of parameters in $f_0$ and LIDs}
  \end{subfigure}
   \begin{subfigure}{.67\linewidth}
    \includegraphics[width=0.82\textwidth]{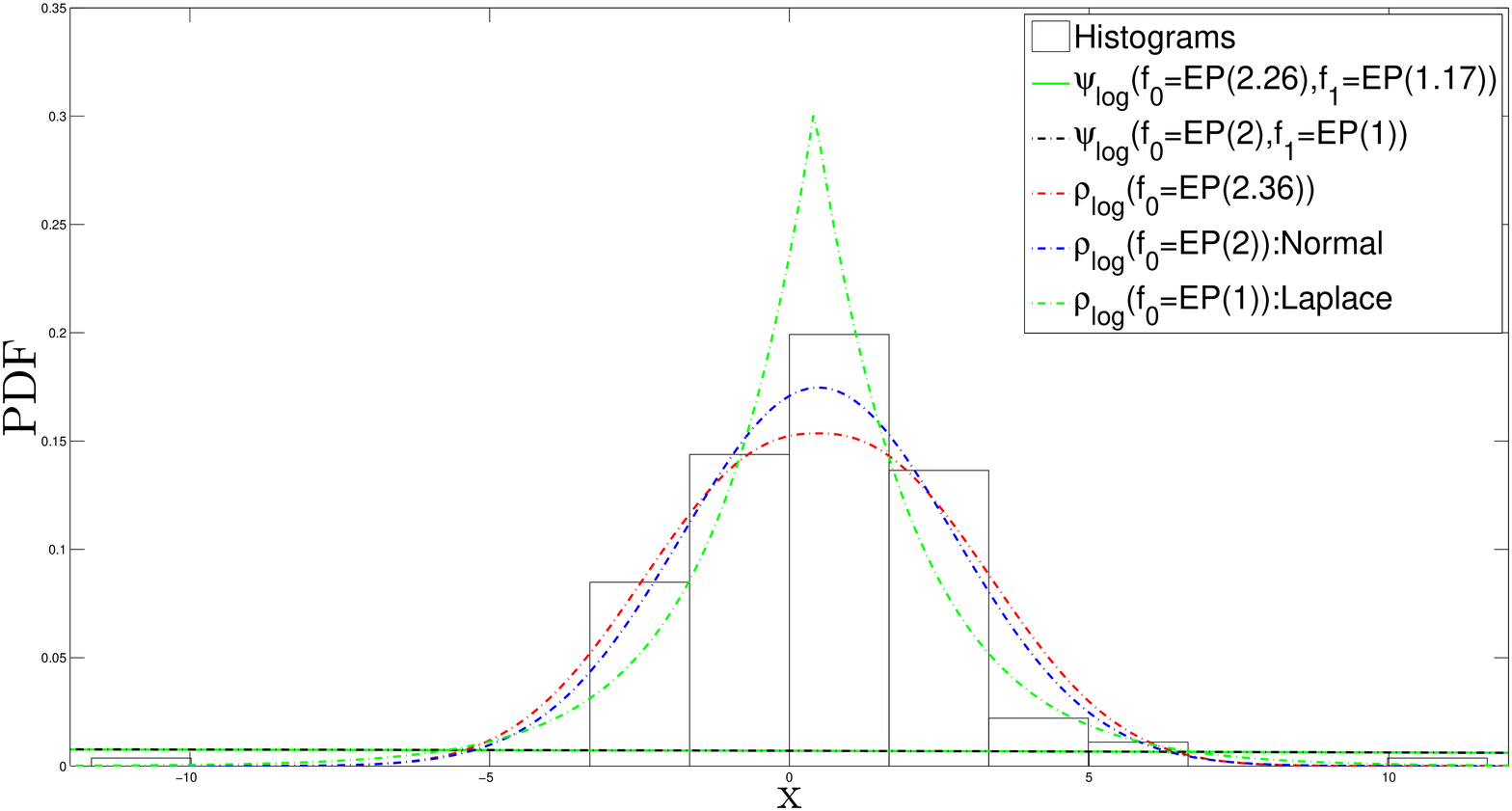}%pdfa2a1collectionDigerleri pdfa2a1collectionDigerleriOUT
    \caption{Non robust M-estimators and MLE of parameters in $f_0$ and LIDs when two outliers are added}
  \end{subfigure}
   \begin{subfigure}{.67\linewidth}
    \includegraphics[width=0.82\textwidth]{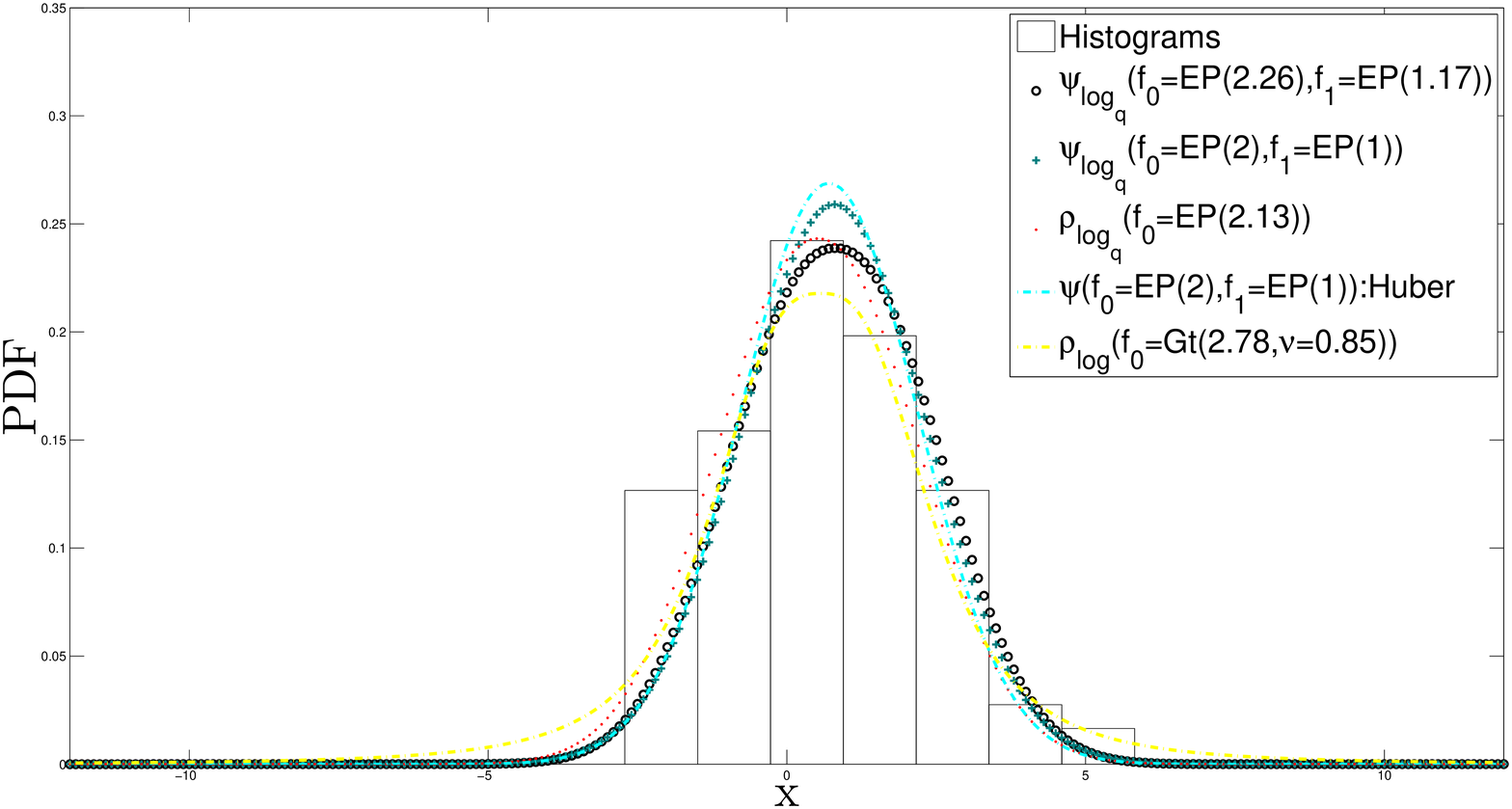}%pdfa2a1collectionDigerleri
    \caption{Robust M-estimators and MLE of parameters in $f_0$ and LIDs}
  \end{subfigure}
   \begin{subfigure}{.67\linewidth}
    \includegraphics[width=0.82\textwidth]{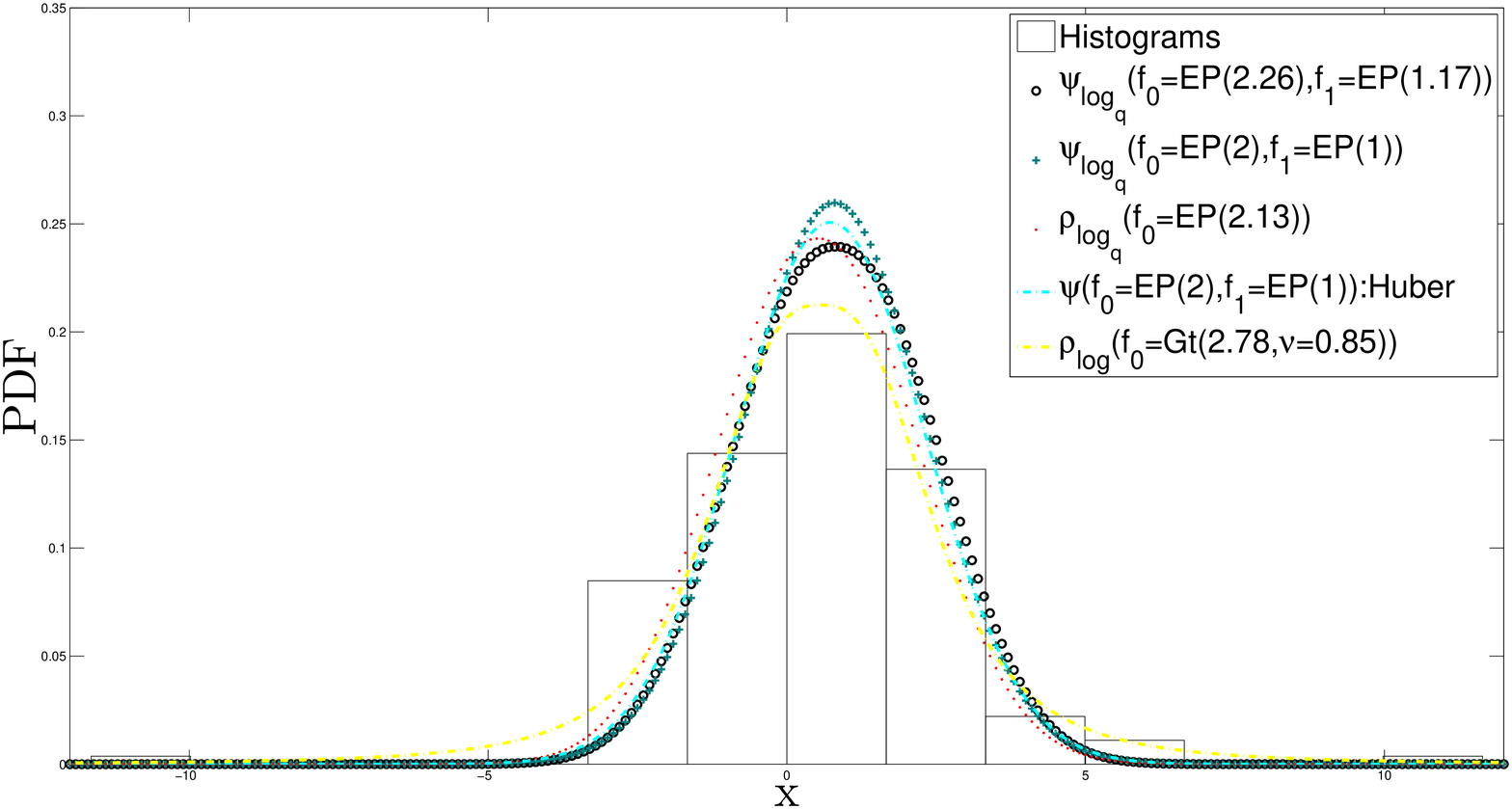}%pdfa2a1collectionDigerleri pdfa2a1collectionDigerleriOUT
    \caption{Robust M-estimators and MLE of parameters in $f_0$ and LIDs when two outliers are added}
  \end{subfigure}
  \caption{M-estimators and MLE of parameters in $f_0$ and LIDs for protein data in cancer cell}
    \label{Ex3figtoplucaMS}
  \end{figure}

Let us think about a comparison between LID and $f_0$. When $\psi_{\log_{q}}(f_0,f_1)$ are compared with $\rho_{\log_{q}}(f_0)$, it is seen that $\rho_{\log_{q}}(f_0)$ of $f_{\text{EP}}$ and $\rho_{\log}(f_0)$ of $f_{\text{Gt}}$ are unimodal distributions, but LID as a combination of the functions $f_0$ and $f_1$ can be a bimodal distribution due to the fact that it represents the mixed distribution. It is observed that $\psi_{\log}(f_0,f_1)$ or $\psi_{\log_q}(f_0,f_1)$ can be preferred in such a situation, because the positive or negative sides of histograms in Figure \ref{Ex3figtoplucaMS} show that there can be $f_1$ in data. Note that the bimodality is an example for non identical distribution. LID as a mixing of $f_0$ and $f_1$ can accomplish to fit data in bimodal case as well.

When we make a comparison between (a) and (b) in Figure \ref{Ex3figtoplucaMS}, $\psi_{\log}(f_0=EP(p_0),f_1=EP(p_1))$ does not  produce robust and efficient estimators. For this reason, LIDs have to be used in $\log_q$.  When we look at (c) and (d) in Figure \ref{Ex3figtoplucaMS}, it is observed that the best fitting can be accomplished by $\psi_{\log_q}(f_0=EP(p_0),f_1=EP(p_1))$, because the data can be modelled well due to $f_1$ in LID if it is compared with $\log_q(f_0)$. It should be also noted that the parameter $p_0$ in $f_0$ and the parameter  $p_1$  in $f_1$ have an important role in getting the efficient M-estimators for the parameters $\mu$ and $\sigma$. Since $f_0$ and $f_1$ are unimodal distributions, $\Lambda(f_0)$ or $\Lambda(f_1)$ is unimodal, because the function $\Lambda$ is an one to one mapping from $f_0$ to $\Lambda$.

  The M-estimators from $\rho_{\log_q}$ are efficient and robust for estimations of $\mu$ and $\sigma$ when they are compared with MLE of $\mu$ and $\sigma$ in exponential power (EP) distribution. Huber M-estimation is sensitive for this data, but $\rho_{\log_q}(f_0)$ are insensitive when their outlier cases are considered in their self, because the parameters $q$, $p_0$
and $p_1$ help us to manage the behaviour of function for efficiency and robustness together. However, the tuning parameter $u$ in Huber M-estimation is not capable of fitting the shape of data set. The parameter $u$  is responsible for determining  where the normal distribution is ended and the Laplace distribution is started. For the information criterions, $RBIC_q$ depends on a part $\log(n)$, therefore $RBIC_q$ can detect whenever there is a change in the sample size of data set.

\begin{table}[!htb]
\centering
\caption{Estimates of parameters $\mu$ and $\sigma$ from different estimating functions without and with two outliers for protein data in cancer cell}
\scalebox{1.05}{
\label{ex1tablemusigma}
\begin{tabular}{ccccc}
 Estimating Functions &$\hat{\mu}$& $\hat{\sigma}$  & $RAIC_{q}$ & $RBIC_{q}$   \\ \hline \hline
$\psi_{\log_{q=.01}}(f_0=EP(2.26),f_1=EP(1.17))$ &  0.8146     &  1.6481  &      20.1088    &  26.2840     \\
Two Outliers&   0.8152  &  1.6434     &  20.0963 &  26.2961    \\\hline
$\psi_{\log_{q=.01}}(f_0=EP(2),f_1=EP(1))$ &   0.7952     &  1.5401  &   22.6947   &  28.8699  \\
Two Outliers   &   0.7956 &   1.5358   &  22.6742 &  28.8740  \\ \hline
$\psi_{\log}(f_0=EP(2.26),f_1=EP(1.17))$&    0.8326     & 4.6304  &  124.5664   & 130.7416    \\
Two Outliers &   -26.8057    & 50.0000  & 134.5075 & 140.7072    \\\hline
$\psi_{\log}(f_0=EP(2),f_1=EP(1))$&    0.8492 & 4.7836 & 140.7803 & 146.9555  \\
Two Outliers &    -24.5185   & 50.0000  &150.4184 & 156.6181  \\ \hline
$\rho_{\log_{q=.66}}(f_0=EP(2.13))$ &  0.5219     &  1.6278  &  470.8822   &  477.0574    \\
Two Outliers&   0.5219 &   1.6277    &    482.6470 &  488.8467  \\\hline
\hline
Estimating Function & $\hat{\mu}$  & $\hat{\sigma}$ & AIC   & BIC    \\ \hline \hline
$\psi(f_0=EP(2),f_1=EP(1),u=1.05)$:Huber&   0.7054     & 1.4840  &   680.5156 & 686.6908  \\
Two Outliers&     0.7308 &   1.5909   & 715.5155 &  721.7153 \\\hline
$\rho_{\log}(f_0=Gt(2.78,\nu=0.85))$ &  0.5661     &  2.2502   &      688.5012 & 694.6764        \\
Two Outliers &     0.5630  &  2.3073   &    717.4041& 723.6038 \\\hline
$\rho_{\log}(f_0=EP(2.36))$&    0.4835    &   2.0091     &   670.6678    &  676.8429  \\
Two Outliers &      0.4703 &   2.5532  & 757.5195  & 763.7192     \\ \hline
$\rho_{\log}(f_0=EP(2))$:Normal&    0.4893       & 1.8978   & 671.3149     & 677.4901         \\
Two Outliers & 0.4833  &  2.2833       &         740.2122 & 746.4119           \\ \hline
$\rho_{\log}(f_0=EP(1))$:Laplace &  0.4232    &  1.5243   & 689.1642 & 695.3394   \\
Two Outliers &       0.4182  &  1.6477  &   723.1529 & 729.3527    \\ \hline
\hline
\end{tabular}}
\end{table}

$\log_q$ is a kind of generalized $t$ as a generalized exponential, that is $\text{exp}_q(x)=[1+(1-q)x^p]^{1/(p-pq)}$ that can be regarded as a similar kernel of Gt in Table \ref{PDFstable} if a kernel $x^p$ of EP is based on $\log_q$ \cite{van95}. The kernel $x^p$ of EP in $\log_q$ has a parameter $q$ that manages the tail behaviour of function as well with overall shape of EP, which shows that the kernel $x^p$ of EP in $\log_q$ can be a good candidate for efficient fitting on data. $\rho_{\log}(f_0=EP(2.36))$ has a good fitting when it is compared with other $\rho_{\log}$ functions, because AIC and BIC of $\rho_{\log}(f_0=EP(2.36))$ are the smallest values among AIC and BIC of other $\rho_{\log}$ functions. From here, it is observed that the nature of data is appropriate to use a function that has flat peakedness  property. In other  words, the data are member of a function, i.e., population, having this property. The flat peakedness can occur due to non-identicality of data as well. In this context, it is observed that LID having parameter $p_0=2.26$ that can model peakedness of data can be preferable. When we consider about fitting performance of $\psi_{\log_{q=.01}}(f_0=EP(2.26),f_1=EP(1.17))$ and $\psi_{\log_{q=.01}}(f_0=EP(2),f_1=EP(1))$, the first one has small AIC and BIC values when it is compared with that of second one for $q=.01$ at a same base for comparison.

For all of four examples, all of Tables show that the estimates of parameters, the values of $RAIC_{q}$ and $RBIC_{q}$ that do not go to big values for cases as the values of AIC and BIC. $RAIC_q$ and $RBIC_q$ from $\psi_{\log_q}(f_0,f_1)$ do not  change much when they are compared with that of $\rho_{\log_q}(f_0)$. MLEs of parameters are not efficient and robust, however M-estimators from $f_0$ and LIDs in  $\log_q$ can be efficient and robust, because we can have neighborhood of PDFs via tuning parameter $q$ and the parameters $p_0$ and $p_1$ in a PDF. Generally, the parameter $q$ in LID case can be near to zero for robustness and efficient fitting. However, $\rho_{\log_q}(f_0)$ case can take values at a large interval for robustness when it is compared with that of $\psi_{\log_q}(f_0,f_1)$ for a such kind of data set in examples in \ref{ex1musig} and \ref{ex2musig}.

$\psi_{\log_q}(f_0,f_1)$  is a different function from $\Lambda(f_0)$ to fit data set. Since we have a mixed function or the convex combination of $f_0$ and $f_1$, i.e, LID, its behaviour on fitting data is not same with $\Lambda(f_0)$, which shows why some results in outlier case from $\psi_{\log_q}(f_0,f_1)$ have different values. Here, it is taken in account that the data and the function have to accommodate each others well, because we need to get more information from data set. In another side, $\psi_{\log_q}(f_0,f_1)$ can be beneficial when we have $f_1$, that is, the data are distributed non-identically. The nonidentical case can also produce a bimodal distributed data, because equation \eqref{eqlnqLID} can be regarded as a mixed distribution. When we look at general results of examples in \ref{ex1musig} and \ref{ex2musig}, it is observed that efficiency and robustness can work together while performing the estimation procedure.

For case of $\mu$ and $\sigma$, generating artificial data is not required, because EP distribution with estimated values of $\mu$ and $\sigma$ and also $p_0=2.26$ from $\psi_{\log_q}(f_0=EP(2.26),f_1=EP(1.17))$ can represent histograms of real data at Figure \ref{Ex3figtoplucaMS}-(c) and (d)  well when it is compared with other PDFs drawn by the estimated values of parameters $\mu$ and $\sigma$ from $\rho_{\log}$ and $\psi$ in (a)-(d) in Figure \ref{Ex3figtoplucaMS}. Using one symmetric distribution to generate data should not be preferable, because there can be an unknown contamination rate $\epsilon$ from $f_1$. For this reason, the bimodality cannot be constructed by means of two mixing distributions when the contamination rate  cannot be known exactly. Note that the artificial data generation from underlying distribution having shape and scale parameters, such as Gamma, Weibull and Burr is easier than the location and scale models when we want to make a cross check between the artificial data and the real data. In location-scale case, we can need to determine the bimodality from the contamination rate. When this discussion has been taken care, it is not possible to make a comparison between artificial data generated from symmetric distribution and real data of phenomena for finite sampling $n$. Therefore, we omit to give the artificial data for two examples in \ref{appBlocsca}. However, we examine the bin ranges at $[-6,-4,-2,-0.5,0.5,2,4,6]$. For such bin ranges, the numbers of observations at bins are   $[0,23,22,40,40,34,3,0]$. According to these numbers, the numbers of observations at tail of negative and positive parts of real line are $45$ and $37$. From here, it is observed that LIDs and Huber M-estimation are capable of overcoming the effect of tail of negative part of real line.  The numbers of negative and positive observations are $62$ and $100$, respectively. Since the number of positive observations is $100$, it is observed that the underlying distribution can be on positive axis. The estimating functions of LID and Huber M-estimation tend to represent the underlying distribution, as it is observed from the estimates of parameters $\mu$ and $\sigma$ from these functions (see Table \ref{ex1tablemusigma}). The importance of  flatness which can occur due to $f_1$ in data has been observed by $\rho(f_0=EP(2.36))$, because the smallest BIC value is given by $\rho_{\log}(f_0=EP(2.36))$. However, in the two outliers case, the BIC is drastically inflated and BIC value of $\rho_{\log}(f_0=EP(2.36))$ is the biggest one (see Table \ref{ex1tablemusigma}).

\subsection{Example 2}\label{ex2musig}

A protein data coded as  ME:UACC-62 from Lysate Array at a website https://discover.nci.nih.gov/cellminer/ is analysed to see tendency (location) and spread (scale) of protein in cancer cell. The maximum value as an outlier is 8.968 at
positive and -8.968 at negative sides of real axis. Therefore, we keep the symmetry of data. In HGA, our prescribed intervals for $\mu$ and $\sigma$ are $[-50,50]$ and $[0,50]$, respectively. The results in here are supported by the results in example  1. Therefore, we omit to rewrite the same comments for example 2. It is also noted that it is expectable to get the results which can be similar framework with example 1, because the examples 1 and 2 can have similar nature in their self.  However, we will give some important results that are included by the results in example 1.

\begin{table}[!htb]
\centering
\caption{Estimates of parameters $\mu$ and $\sigma$ from different estimating functions without and with two outliers for protein data in cancer cell}
\scalebox{1.05}{
\label{ex2tablemusigma}
\begin{tabular}{ccccc}
 Estimating Functions &$\hat{\mu}$& $\hat{\sigma}$  & $RAIC_{q}$ & $RBIC_{q}$   \\ \hline \hline
$\psi_{\log_{q=.005}}(f_0=EP(2.14),f_1=EP(1.09))$ &  0.5566    &  1.9963  &       19.0393 &  25.2145       \\
Two Outliers&       0.5586  &  1.9863    &  18.9891 &   25.1889       \\\hline
$\psi_{\log_{q=.005}}(f_0=EP(2),f_1=EP(1))$ &   0.5784     &   1.8387  &  20.3508 &  26.5260  \\
Two Outliers   &    0.5795 &   1.8291    &  20.2997 &  26.4994  \\ \hline
$\psi_{\log}(f_0=EP(2.14),f_1=EP(1.09))$&    -2.0707   & 7.0343  &    133.6396 & 139.8148      \\
Two Outliers &      26.5503     & 50.0000 &   141.6468 &  147.8465     \\\hline
$\psi_{\log}(f_0=EP(2),f_1=EP(1))$&  -1.9874    &  7.0704 & 142.6084 &  148.7836  \\
Two Outliers &    25.2716     & 50.0000  & 150.4277  & 156.6274  \\ \hline
$\rho_{\log_{q=.51}}(f_0=EP(2.43))$ &  0.3452   &  1.7677 &    420.3556 &  426.5308    \\
Two Outliers&       0.3453  &    1.7677    &     428.5159  & 434.7156  \\\hline
\hline
Estimating Function & $\hat{\mu}$  & $\hat{\sigma}$ & AIC   & BIC    \\ \hline \hline
$\psi(f_0=EP(2),f_1=EP(1),u=1.08)$:Huber&      0.4722  & 1.5653  &    708.2074 & 714.3826  \\
Two Outliers&   0.5160 &   1.6890   &   735.6557 & 741.8554 \\\hline
$\rho_{\log}(f_0=Gt(2.34,\nu=1.75))$ &   0.3079    &  2.6362   &       704.3447 &  710.5199       \\
Two Outliers &       0.3048  &  2.7208   &     731.4589  & 737.6586 \\\hline
$\rho_{\log}(f_0=EP(2.65))$&       0.2330  &  2.1809    &     684.4131 & 690.5883    \\
Two Outliers &    0.2239  &  2.5449    & 743.4512  & 749.6509      \\ \hline
$\rho_{\log}(f_0=EP(2))$:Normal&       0.2756  &  2.0097   &   689.8820 & 696.0572         \\
Two Outliers &    0.2723 &   2.2296      &      732.4160 &  738.6157     \\ \hline
$\rho_{\log}(f_0=EP(1))$:Laplace &     0.3755 &   1.6327   & 711.4166 & 717.5918   \\
Two Outliers &         0.3989 &   1.7222    &    737.6460  & 743.8457   \\ \hline
\hline
\end{tabular}}
\end{table}

\begin{figure}[!]
\centering
   \begin{subfigure}{.67\linewidth}
    \includegraphics[width=0.82\textwidth]{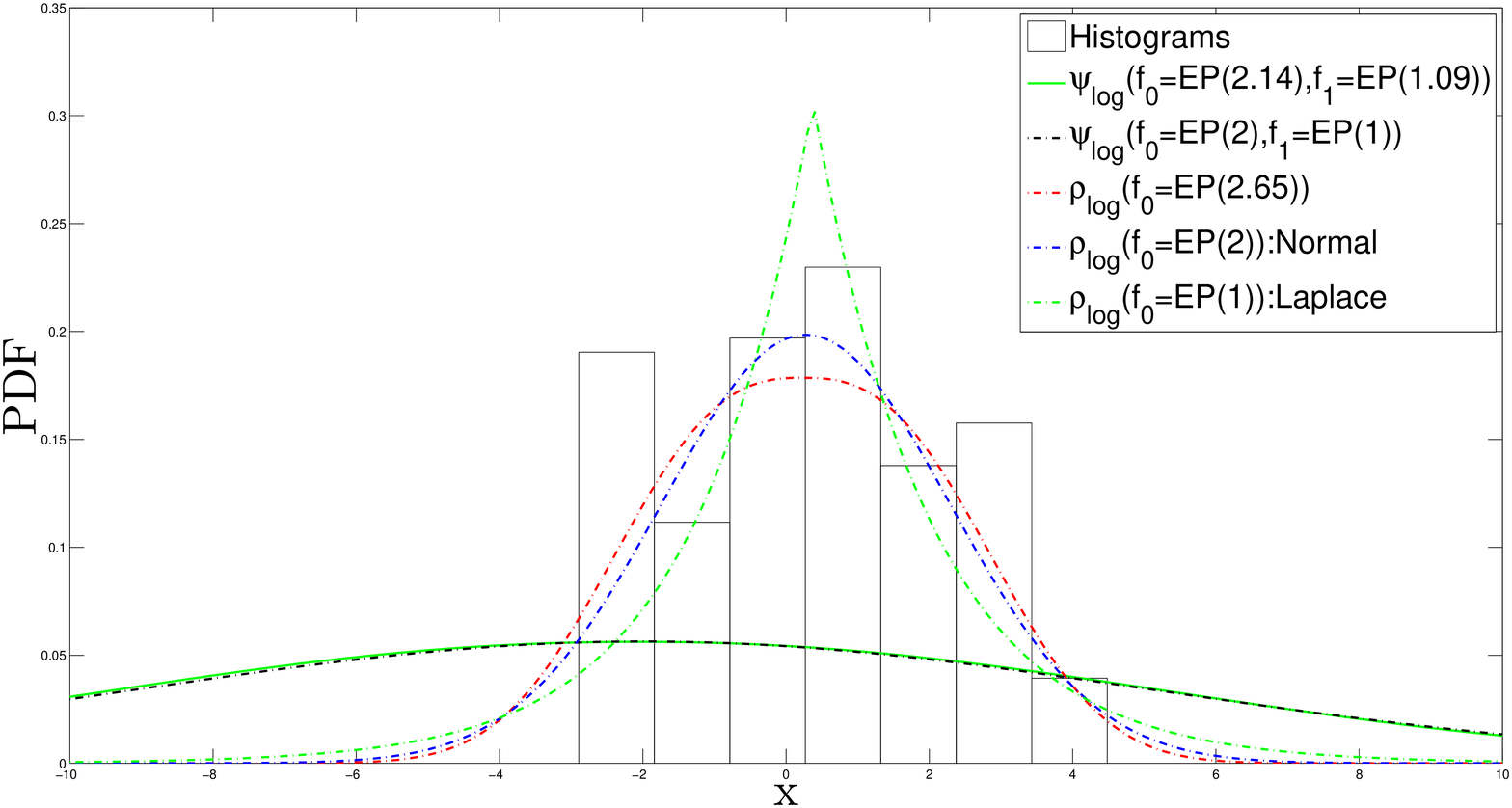}%pdfa2a1collectionDigerleri
    \caption{Non robust M-estimators and MLE of parameters in $f_0$ and LIDs}
  \end{subfigure}
   \begin{subfigure}{.67\linewidth}
    \includegraphics[width=0.82\textwidth]{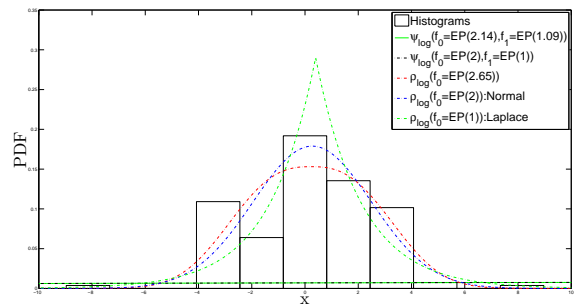}%pdfa2a1collectionDigerleri pdfa2a1collectionDigerleriOUT
    \caption{Non robust M-estimators and MLE of parameters in $f_0$ and LIDs when two outliers are added}
  \end{subfigure}
   \begin{subfigure}{.67\linewidth}
    \includegraphics[width=0.82\textwidth]{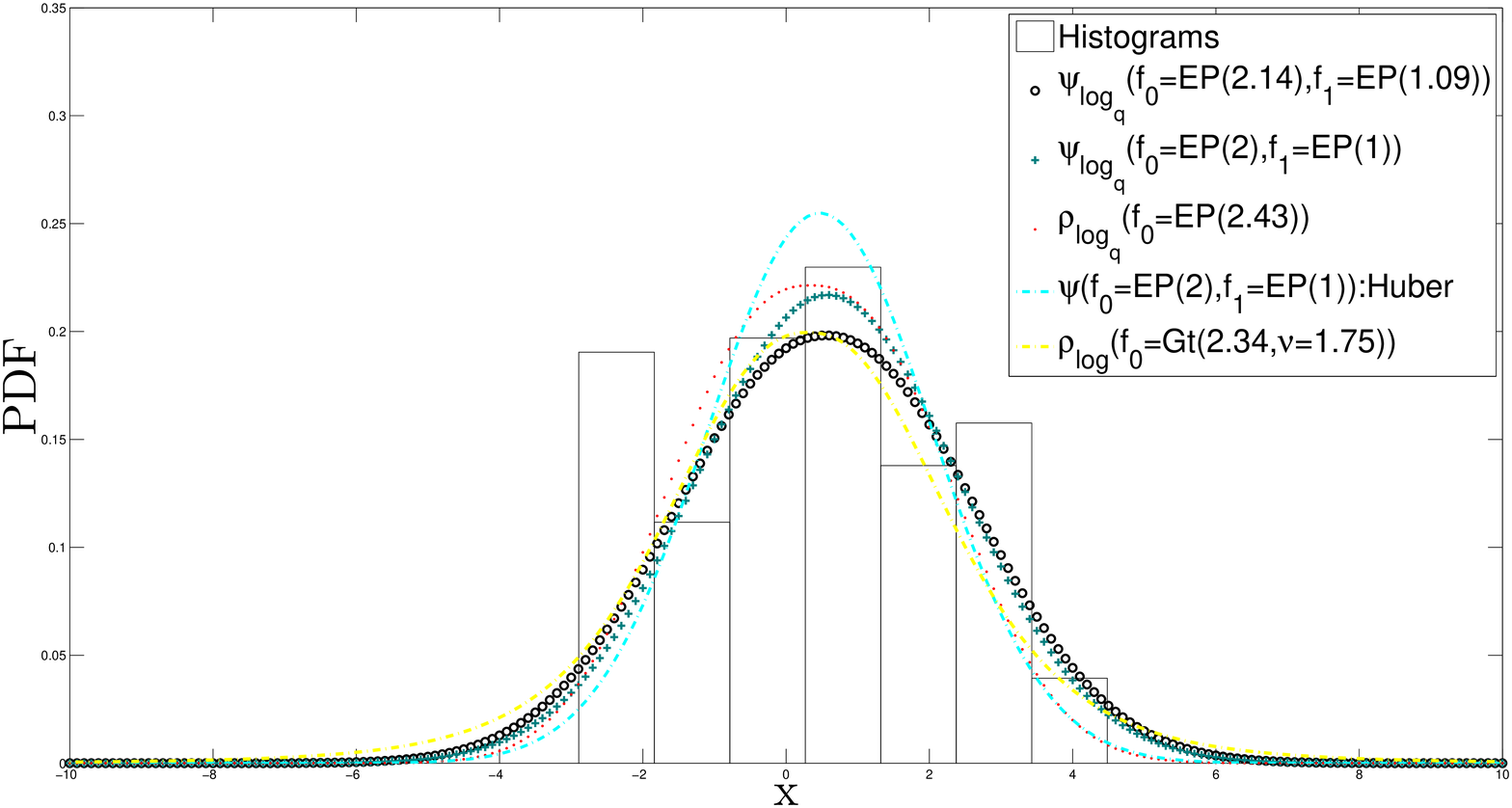}%pdfa2a1collectionDigerleri
    \caption{Robust M-estimators and MLE of parameters in $f_0$ and LIDs}
  \end{subfigure}
   \begin{subfigure}{.67\linewidth}
    \includegraphics[width=0.82\textwidth]{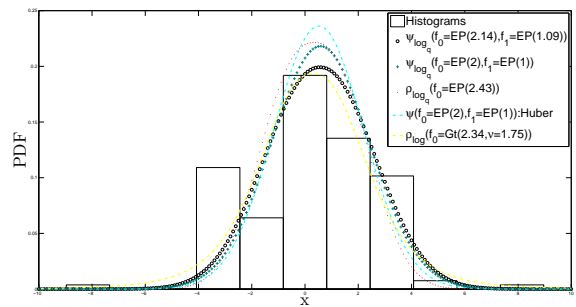}%pdfa2a1collectionDigerleri pdfa2a1collectionDigerleriOUT
    \caption{Robust M-estimators and MLE of parameters in $f_0$ and LIDs when two outliers are added}
  \end{subfigure}
  \caption{M-estimators and MLE of parameters in $f_0$ and LIDs for protein data in cancer cell}
    \label{Ex4figtoplucaMS}
  \end{figure}

We examine the bin ranges at $[-6,-4,-2,-0.5,0.5,2,4,6]$. For such bin ranges, the numbers of observations at bins are   $[0,29,24,33,41,33,2,0]$. According to these numbers, the numbers of observations at tail of negative and positive parts of real line are $53$ and $35$. From here, it is observed that LIDs and Huber M-estimation are capable of overcoming the effect of tail of negative part of real line.  The numbers of negative and positive observations are $69$ and $93$, respectively. Since the number of positive observations is $93$, it is observed that the underlying distribution can be on positive axis. The estimating functions of LID and Huber M-estimation tend to represent the underlying distribution, as it is observed from the estimates of parameters $\mu$ and $\sigma$ from these functions. The importance of  flatness which can occur due to $f_1$ in data has been observed by $\rho(f_0=EP(2.65))$, because the smallest BIC value is given by $\rho_{\log}(f_0=EP(2.65))$. However, in the two outliers case, the BIC is drastically inflated and BIC value of $\rho_{\log}(f_0=EP(2.65))$ is the biggest one (see Table \ref{ex2tablemusigma}).

Note that the non identically distributed data, namely $f_1$ in data, will affect the estimations of parameters. Thus, LID is beneficial for overcoming the problem of modelling $f_1$ as well and the importance of LID can be observed. Here, $q$ and $p$ in LID case (or $\nu$ and $p$ in Gt distribution) can interact with each other as shape parameters. Therefore, we can also need LID case for modelling data efficiently via the functions $f_0$ and $f_1$ to render the effect of this interaction as possible as we can do. In fact, the results of estimates of $\mu$ support that LID and Huber with their functions $f_0$ and $f_1$ are better than Gt with only $f_0$ when the modelling performance for data set which can be considered to come from the underlying distribution $f_0$ and the contamination with $f_1$ is taken in account. The estimates of location parameter tend to go to the positive side of data which can be considered to be a member of the $f_0$.  

Especially, the estimation process requires to have LID, because we have data, that is, there is a sampling version of unknown function. Since we handle with unknownness, using functions that can be neighborhood to each others via q and LID will help us to have precise and robust estimated values of parameters, as it is observed from the applications on real data sets. At the end, the position of data where they are and the position of function where it is are extremely important to have efficient M-estimators for parameters $ \boldsymbol \theta $.

\end{document}